\documentclass{article}
\usepackage[margin=1in]{geometry}
\usepackage{microtype}
\usepackage{amsmath,amssymb,amsthm}
\newtheorem{theorem}{Theorem}
\newtheorem{lemma}{Lemma}

\usepackage[square,numbers,sort&compress,longnamesfirst]{natbib}
\usepackage{hyperref}
\usepackage[nameinlink,noabbrev]{cleveref}
\Crefformat{equation}{(#2#1#3)}
\labelcrefformat{equation}{(#2#1#3)}
\Crefrangeformat{equation}{(#3#1#4)--(#5#2#6)}
\Crefrangeformat{equation}{(#3#1#4)--(#5#2#6)}
\Crefmultiformat{equation}{(#2#1#3)}{ and~(#2#1#3)}{, (#2#1#3)}%
  {, and~(#2#1#3)}
\Crefmultiformat{equation}{(#2#1#3)}{ and~(#2#1#3)}{, (#2#1#3)}%
  {, and~(#2#1#3)}

\usepackage{xspace}
\newcommand{\ie}{i.e.,\xspace}
\newcommand{\eg}{e.g.,\xspace}
\newcommand{\E}{\mathbb{E}}

\newcommand{\Prob}{\mathbb{P}}
\newcommand{\ind}{\mathbf{1}}
\newcommand{\dd}{\,\mathrm{d}}

\newcommand{\calT}{\mathcal{T}}
\newcommand{\calX}{\mathcal{X}}
\newcommand{\sbar}{\bar{s}}

\newcommand{\BS}{\mathrm{BS}}
\newcommand{\TS}{\mathrm{TS}}
\newcommand{\SEC}{\mathrm{S}}

\usepackage[dvipsnames]{xcolor}
\usepackage{todonotes}

\title{Prophet Inequalities for Mixtures of I.I.D.\@ Random Variables}
\author{Felix Fischer\thanks{School of Mathematical Sciences, Queen Mary University of London, UK, \href{mailto:felix.fischer@qmul.ac.uk}{felix.fischer@qmul.ac.uk}.} \and Kevin Schewior\thanks{Department of Mathematics and Computer Science, University of Cologne, Germany, \href{mailto:k.schewior@uni-koeln.de}{k.schewior@uni-koeln.de}.}}
\date{}

\begin{document}

\maketitle

\begin{abstract}
A latent index $J\in\{1,\ldots,k\}$ is drawn from a prior~$\pi$. Conditional on $J=j$, a gambler observes random variables $X_1,\ldots,X_n$ drawn independently from a non-negative component distribution $F_j$. The gambler knows $n$, the prior $\pi$, and the component distributions $F_1,\ldots,F_k$, but not $J$, observes the values sequentially, and selects at most one value at the time it is revealed.
We show that the largest $\alpha$ such that the gambler can guarantee to
select a value of $\alpha\cdot\mathbb{E}[\max_iX_i]$ in expectation, for all
distributions where this expectation is finite, is
\[
  \alpha=\frac1e+\frac{1}{k^{2/3+o(1)}}.
\]
The lower bound is achieved by the better of two variants of the standard secretary rule. The upper bound uses a finely tuned instance of the common-scale model, which scales a sequence of independent random variables by a single latent variable.
\end{abstract}

\section{Introduction}

For a given sequence of random variables, a prophet inequality compares the value achievable by a ``gambler,'' who observes the sequence in order and selects an element at the time it is revealed to the value achievable by a ``prophet,'' who knows the entire sequence~\citep[\eg][]{prophetbook}. We consider the case where the sequence is a mixture of i.i.d.\@ draws, and elements of the sequence are correlated from the perspective of the gambler due to lack of knowledge of the particular component of the mixture they are drawn from.

Let $J\in\{1,\ldots,k\}$ be drawn from a prior $\pi$. Conditional on $J=j$, the gambler observes random variables $X_1,\ldots,X_n$ drawn independently from a non-negative component distribution $F_j$. The gambler knows $n$, the prior $\pi$, and the component distributions $F_1,\ldots,F_k$, but not the value of~$J$, observes the values sequentially, and selects at most one value at the time it is revealed. Setting $X_\infty=0$, the gambler's strategy can be described by a random variable $\tau\in\{1,\dots,n\}\cup\{\infty\}$, the so-called stopping time, such that the event $\{\tau=i\}$ for $i=1,\dots,n$ is determined by the observations $X_1,\dots,X_i$ and any private randomness of the gambler. Let $M=\max_{i\le n}X_i$ be the value achievable by the prophet, $V=\sup_\tau\E[X_\tau]$ the value achievable by the gambler, where the expectation is taken over the randomness of the draws and possible randomness of the stopping time. We will be interested in determining the \emph{competitive ratio} $c_k=\inf_{n,\pi,F_1,\ldots,F_k:0<\E[M]<\infty}\frac{V}{\E[M]}$ achieved by the best strategy of the gambler, as a function of the number~$k$ of component distributions.

This problem was first considered by \citet{CorreaDFSZ2021} and interpolates between two settings that are fully understood. Setting $k=1$ recovers the prophet problem for a sequence of i.i.d.\@ random variables~\citep{HillKertz1982}, and thus $c_1\approx 0.745$~\citep{CorreaFHOV2021}. In the absence of a bound on $k$, the best guarantee that can be achieved is $1/e$~\citep{CorreaDFSZ2021}, which implies $\lim_{k\to\infty}c_k=1/e$. We obtain the following result on the asymptotic behavior of $c_k$.
\begin{theorem}
\label{thm:main}
There exist absolute constants $c,C>0$ such that, for all $k\geq 2$, %
\begin{equation*}
 \frac{1}{e}+\frac{c}{(k\log k)^{2/3}}
 \leq c_k \leq
 \frac{1}{e}+\frac{C}{k^{2/3}}.
\end{equation*}
In particular $c_k\to 1/e$ and $c_k-1/e=k^{-2/3+o(1)}$.
\end{theorem}

We will prove the lower bound in \Cref{sec:lower} and the upper bound in \Cref{sec:upper}.

\citet{AmaniHamedaniAIKS2026} have, concurrently and independently, shown a lower bound of $1/e + d/k^2$ and an upper bound of $1/e+D/\sqrt{k}$ for absolute constants $d$ and $D$. They also consider more general problems where random variables need not be identically distributed or independent. Their lower bound has the additional property that it can be approximated by a stopping rule that is computationally efficient; we make no such claim about our lower bound.

\section{Lower Bound}
\label{sec:lower}

We first recall that it is without loss to assume that distributions are continuous. This can be seen by perturbing the values by independent uniform values in $[0,\eta]$ for some $\eta>0$, working with the resulting continuous distribution, and letting $\eta$ go to zero~\citep{prophetbook}. For a CDF $F$ and $0<u<1$ let $F^{-1}(u)=\inf\{x\in\mathbb R:F(x)\ge u\}$, and observe that by continuity $F(F^{-1}(u))=u$.

The lower bound will be achieved by using the better of two stopping rules on any given mixture depending on a parameter of that mixture. Both rules are variants of the optimal rule for the secretary problem; the secretary rule rejects an (expected) initial $1/e$ fraction of the $n$ values and thereafter stops at the first \emph{record}, \ie the first value that exceeds all previous values~\citep{Ferguson1989}. Both variants cut the supports of distributions in the mixture into a set of intervals but then use these cuts in different ways.

Fix $L\geq 1$. For $j=1,\dots,k$ and $\ell=0,\dots,L$, let $T_{j,\ell}=F_j^{-1}(1-1/(2^{\ell}n))$, and observe that by continuity $\Prob_{X\sim F_j}(X>T_{j,\ell})=1/(2^{\ell}n)$. Let $\calT=\{T_{j,\ell}:1\leq j\leq k, 0\leq\ell\leq L\}$. Then $\calT$ splits the support of the distributions into $N+1$ intervals, or \emph{buckets}, where $N\leq k(L+1)$ (with equality if all $T_{j,\ell}$ are distinct).
The first variant, the \emph{bucket secretary rule}, rejects an initial $1/e$ fraction of the values and then stops at the first \emph{bucket record}, \ie the first value occupying a higher bucket than all previous values.
The second variant, the \emph{threshold secretary rule}, uses a threshold $T$ chosen uniformly at random from $\calT$, and a second deadline $b\leq 1/e$; it rejects an initial~$b$ fraction of the values, then accepts the first record above~$T$ until it has seen a $1/e$ fraction of the values, and thereafter accepts the first record.
The bucket secretary rule has to our knowledge not previously been studied. Other authors have combined the secretary rule with a threshold applied to an initial portion of the samples, specifically one that comes from data rather than a distribution~\citep{AntoniadisGouleakisKleerKolev2020,FengLiLiWuWu2025}. Our choice of threshold, sampled from quantile cuts of the component distributions, and our analysis, relating its benefit to the same bucket-gap parameter that controls the bucket secretary rule, are new.

Let
\[
 a=1/e,\qquad P=\E[M],\qquad \gamma=\frac{\E[M_2]}{P},
\]
where $M_2$ is the second-largest of the values $X_1,\dots,X_n$. It is without loss to assume that the prior assigns positive probability to every component, and thus that every component has a finite expected maximum. For $n\le 2$, choosing a uniform observation achieves competitive ratio at least $1/2$. For $n\ge 3$ the guarantee we obtain will be the maximum of the competitive ratios of the two stopping rules. 
The secretary rule's uniform guarantee of~$a$ will serve as a baseline, but we will see that one of the two variants improves on it for every mixture. We review the analysis of the secretary rule in \Cref{sec:secretary}.
Which of the variants performs better will depend on a parameter~$G$ of the mixture, which we will define formally later. The parameter measures the gap, in terms of the number of buckets and weighted by the prophet's reward $M$, between the highest and second-highest occupied buckets. The geometric sizes of the buckets mean that the highest occupied bucket will usually be singly occupied, so that bucketing only incurs a small loss; we will show this in \Cref{sec:loss}. When $G$ is small, occupied buckets below the highest occupied one will often contain multiple observations; seeing one of them early will prevent the bucket secretary rule from stopping too soon. When $G$ is large, a random cut will often separate the maximum value from all other values, allowing the threshold secretary rule to select it before time~$1/e$. We will analyze the bucket secretary rule in \Cref{sec:bucket}, the threshold secretary rule in \Cref{sec:threshold}, and combine the two to obtain the lower bound in \Cref{sec:gain}.

\subsection{Arrival Times and the Secretary Rule}
\label{sec:secretary}

We will use virtual arrival times to analyze stopping rules uniformly over all finite $n$. Independently of the observations, generate $n$ independent uniform random variables on $[0,1]$, sort them, and assign the $i$th smallest time to the $i$th observation. Conditional on the unordered sample, exchangeability makes the observation order uniform. Thus the resulting timed sequence has the same distribution as one obtained by assigning independent uniform arrival times to the labeled observations in the sample and revealing them chronologically. The gambler can implement this construction using private randomness, and we henceforth describe stopping rules in terms of these virtual arrival times.

Let $\SEC_t$ be the stopping rule that rejects before time $t\in[0,1]$ and then accepts the first record. Fix an unordered sample of distinct values and assign independent uniform arrival times to the values. Conditional on the maximum value in the sequence arriving at $u\ge t$, the largest observation to arrive before time $u$, if there is such an observation, has a uniform arrival time in $[0,u]$. $\SEC_t$ thus selects the maximum if that time is below $t$, and with certainty if there is no earlier observation, and therefore
\[
 \Prob(\text{$\SEC_t$ selects $M$})\ge \int_t^1\frac{t}{u}\dd u=-t\log t.
\]
This lower bound is maximized at $t=a=1/e$. We will call $\SEC_a$ the secretary rule; it is the timed version of the classical asymptotically optimal rule for the secretary problem and guarantees selection of the maximum with probability at least $1/e$ for every finite $n$.
Similarly, let $0<t\le a$ and condition on the second maximum
arriving at time $u\ge a$. The maximum arrives later with probability $1-u$, and
the same reasoning applied to the remaining $n-2$ observations yields
\begin{equation}\label{eq:second-max}
 \Prob\bigl(\text{$\SEC_t$ selects the second maximum at a time~$\ge a$}\bigr)
 \ge t\int_a^1\frac{1-u}{u}\dd u=at.
\end{equation}
Both bounds hold for any sample of distinct values, and selection of the maximum and the second maximum are disjoint events. It follows, for example, that $\E[\SEC_a]/P\ge a+a^2\gamma$.

\subsection{Loss Through Bucketing Is Small}
\label{sec:loss}

Let $\delta=2^{-L}$, and recall $N=|\mathcal{T}|\leq k(L+1)$ is the number of cuts of the value space into buckets. Index the buckets $0,\ldots,N$. Call a sample \emph{good} if its highest occupied bucket contains exactly one observation, and let $q$ be the fraction of $M$ coming from good samples, \ie
\[
 q=\frac{\E[M\ind_{\mathrm{good}}]}{P}.
\]
Good samples are good for the bucket secretary rule, because the maximum value is the only value in its bucket and thus remains available for selection. The following lemma bounds the fraction of the prophet’s reward contributed by bad samples in terms of $\gamma$ and $\delta$.
\begin{lemma}
\label{lem:loss}
The share of the prophet's reward on bad samples satisfies
\begin{equation}\label{eq:q}
 1-q\le\gamma+\delta.
\end{equation}
\end{lemma}
\begin{proof}
Fix a component~$j$. Define the quantiles $U_i=F_j(X_i)$ and $U_{\max}=\max_iU_i$, and observe that $U_1,\dots,U_n$ are independent uniform random variables on $[0,1]$ and $M=F_j^{-1}(U_{\max})$ almost surely. Condition on $U_{\max}=u$, and observe that the remaining $n-1$ quantiles are independent uniform draws from $[0,u]$. We will bound the collision probability that the maximum is \emph{not} alone in its bucket and then integrate over $u$.

If $u\le1-1/n$, then $(n-1)(1-u)/u\geq 1$. If $1-1/n<u\le1-\delta/n$, the geometric spacing ensures that the part of the maximum's bucket below $u$ has mass at most $1-u$. The collision probability is therefore at most $(n-1)(1-u)/u$ in both cases, trivially in the former and by the union bound in the latter. If $u>1-\delta/n$, the mass is at most $\delta/n$, so the collision probability is at most $(n-1)\delta/(nu)\le\delta$.
Since $U_{\max}$ has density $nu^{n-1}$ and the maximum reward is $F_j^{-1}(u)$,
\begin{align*}
 \E[M\ind_{\rm bad}\mid J=j]
 &\le\int_0^1F_j^{-1}(u)n(n-1)u^{n-2}(1-u)\,du+\delta\,\E[M\mid J=j] \\
 &=\E[M_2\mid J=j]+\delta\,\E[M\mid J=j],
\end{align*}
where the last equality uses that, conditioned on $j$, the second-largest quantile has density $n(n-1)u^{n-2}(1-u)$, the density of $\mathrm{Beta}(n-1,2)$. Taking expectation over the prior yields
\[
 \E[M\mathbf 1_{\mathrm{bad}}]
\le \E[M_2]+\delta\,\E[M]
\]
and dividing by $P=\E[M]$ shows the claim.
\end{proof}

\subsection{Small Gap Helps Bucket Secretary}
\label{sec:bucket}

For a good sample let $D$ denote the difference between the indices of the largest and second-largest occupied buckets, and let
\begin{equation*}
 G=\frac{\E[MD\ind_{\mathrm{good}}]}{P}.
\end{equation*}
$G$ is a parameter of the instance which measures, weighted by the prophet’s reward, the expected gap between the top two occupied buckets. The bucket secretary rule $\BS$ rejects before time $a$ and then accepts the first observation that sits in a higher bucket than all previous observations, including rejected ones. When the gap $G$ is small, occupied buckets below the highest occupied one tend to contain multiple observations; seeing one of them early prevents $\BS$ from stopping too soon.

On good samples, $1\le D\le N$ and thus $q\le G\le Nq$. In particular $G>0$ whenever $q>0$. On a good sample where the highest occupied bucket below the maximum contains at least two observations, $\BS$ selects the maximum with an additional probability of $\kappa=a(1-a)^2/2$ compared to the secretary rule. This increases the contribution from good samples. The following lemma balances this gain against the fraction of the prophet’s reward contributed by bad samples.
\begin{lemma}
\label{lem:bucket}
If $q>0$, then
\begin{equation*}
 \frac{\E[\BS]}{P}-a\ge\frac{\kappa q^2}{4G}-a(1-q).
\end{equation*}
\end{lemma}
\begin{proof}
Let $R$ be the event that the highest occupied bucket contains exactly one observation and the second-highest contains at least two observations. For $j\in\{1,\ldots,k\}$ and $m\in\{1,\ldots,N\}$, let $E_{j,m}$ be the event that $J=j$, exactly one observation lies in bucket $m$, and the remaining observations lie in lower buckets. These events partition the good samples.

Fix an event $E=E_{j,m}$ of positive probability. Conditional on the identity of the unique observation in bucket $m$, its value $M$ has distribution $F_j$ restricted to bucket $m$, while the remaining observations are independent draws from $F_j$ restricted to buckets below $m$, independent of $M$. Label these remaining observations in their original index order. Their joint distribution with $M$ does not depend on the identity of the observation in bucket $m$, so the same independence holds conditional on $E$ alone. Let $Z$ be the index of their highest occupied bucket, which is the highest occupied bucket below $m$, and let $p_E=\Prob(R\mid E)$ be the probability that bucket $Z$ contains more than one observation.

We will first prove a lower bound on $p_E$. The $n-1$ observations below the maximum are independent draws from $F_j$ conditioned to lie below bucket $m$, so their bucket indices are drawn independently from some distribution $H$ on $\{0,\ldots,m-1\}$. Let $h_\ell=H(\ell)^{(n-1)/2}=\sqrt{\Prob(Z\leq\ell\mid E)}$ and observe that $h_{-1}=0$.
For each $\ell\in\{0,\ldots,m-1\}$, by independence,
\[
 \Prob(Z=\ell\mid E) = \Prob(Z\le\ell\mid E)-\Prob(Z\le\ell-1\mid E) = h_\ell^2-h_{\ell-1}^2.
\]
Since there are $n-1$ choices for the largest value, it lies in bucket~$\ell$ with probability $H(\ell)-H(\ell-1)$, and each of the remaining $n-2$ observations lies in a lower bucket with probability $H(\ell-1)$,
\begin{align*}
 &\Prob(Z=\ell, \text{exactly one observation in bucket $\ell$}\mid E)\\
 &\qquad=(n-1)\bigl(H(\ell)-H(\ell-1)\bigr)H(\ell-1)^{n-2}.
\end{align*}
When $H(\ell-1)>0$, convexity of $x^{(n-1)/2}$ for $x\geq 0$ yields
\begin{align*}
 h_\ell-h_{\ell-1}
 &=H(\ell)^{(n-1)/2}-H(\ell-1)^{(n-1)/2}\\
 &\ge\frac{n-1}{2}\bigl(H(\ell)-H(\ell-1)\bigr)
 H(\ell-1)^{(n-3)/2}
\end{align*}
and thus
\begin{equation}
 \label{eq:lower}
\begin{aligned}
 2h_{\ell-1}(h_\ell-h_{\ell-1})
 &\ge(n-1)\bigl(H(\ell)-H(\ell-1)\bigr)
 H(\ell-1)^{(n-1)/2+(n-3)/2}\\
 &=(n-1)\bigl(H(\ell)-H(\ell-1)\bigr)H(\ell-1)^{n-2}.
\end{aligned}
\end{equation}
When $H(\ell-1)=0$,~\Cref{eq:lower} holds because both sides are zero. It follows that
\begin{align*}
 &\Prob(Z=\ell,\ \text{at least two observations in bucket $\ell$}\mid E)\\
 &\quad=\Prob(Z=\ell\mid E)
 -\Prob(Z=\ell,\ \text{exactly one observation in bucket $\ell$}\mid E)\\
 &\quad=h_\ell^2-h_{\ell-1}^2
 -(n-1)\bigl(H(\ell)-H(\ell-1)\bigr)H(\ell-1)^{n-2}\\
 &\quad\ge h_\ell^2-h_{\ell-1}^2-2h_{\ell-1}(h_\ell-h_{\ell-1})\\
 &\quad=h_\ell^2-2h_\ell h_{\ell-1}+h_{\ell-1}^2\\
 &\quad=(h_\ell-h_{\ell-1})^2,
\end{align*}
where the inequality holds by~\Cref{eq:lower}, and thus
\begin{align}
 \label{eq:plower}
 p_E=\sum_{\ell=0}^{m-1}
 \Prob(Z=\ell,\ \text{at least two observations in bucket $\ell$}\mid E)
 \ge\sum_{\ell=0}^{m-1}(h_\ell-h_{\ell-1})^2.
\end{align}
Also,
\begin{align}
 \label{eq:mminus}
 \sum_{\ell=0}^{m-1}h_\ell^2
 =\sum_{\ell=0}^{m-1}\Prob(Z\le\ell\mid E)
 =\E\left[\sum_{\ell=0}^{m-1}\ind_{\{Z\le\ell\}}\mid E\right] 
 =\E[m-Z\mid E]=m-\E[Z\mid E].
\end{align}
Therefore
\begin{align*}
 1&=h_{m-1}^2-h_{-1}^2\\
 &=\sum_{\ell=0}^{m-1}(h_\ell^2-h_{\ell-1}^2)\\
 &=\sum_{\ell=0}^{m-1}(h_\ell-h_{\ell-1})(h_\ell+h_{\ell-1})\\
 &\le\left(\sum_{\ell=0}^{m-1}(h_\ell-h_{\ell-1})^2\right)^{1/2}
      \left(\sum_{\ell=0}^{m-1}(h_\ell+h_{\ell-1})^2\right)^{1/2}\\
 &\le\sqrt{p_E}\left(4\sum_{\ell=0}^{m-1}h_\ell^2\right)^{1/2}\\
 &=2\sqrt{p_E(m-\E[Z\mid E])},
\end{align*}
where the first equality holds because $h_{m-1}=1$ and $h_{-1}=0$, the second equality by telescoping, the first inequality by Cauchy--Schwarz, the second inequality by~\Cref{eq:plower} and because $0\le h_{\ell-1}\le h_\ell$, and the last equality by~\Cref{eq:mminus}.
It follows that
\[
 p_E\ge\frac{1}{4(m-\E[Z\mid E])}.
\]
Define
\[
d_E=\E[D\mid E]=m-\E[Z\mid E],
\]
so that $p_E\ge 1/(4d_E)$.

We will now average over all positive-probability events $E=E_{j,m}$ to obtain a lower bound on the fraction of the prophet’s expected reward contributed by samples in $R$.
Let
\[
 w_E=\frac{\Prob(E)\E[M\mid E]}P.
\]
Since the events partition the good samples,
\begin{align*}
 \sum_Ew_E
 &=\frac1P\sum_E\Prob(E)\E[M\mid E]\\
 &=\frac{\E[M\ind_{\rm good}]}P=q.
\end{align*}
For each $E$, $M$ is independent of the other observations, while $D$ and $R$ depend only on their bucket indices. Thus
\begin{align*}
 \sum_Ew_Ed_E
 &=\frac{1}{P}\sum_E\Prob(E)\;\E[M\mid E]\;\E[D\mid E] \\
 &=\frac{1}{P}\sum_E\Prob(E)\;\E[MD\mid E] \\
 &=\frac{\E[MD\ind_{\rm good}]}{P}=G,
\end{align*}
and
\begin{align*}
 \frac{\E[M\ind_R]}P
 &=\frac1P\sum_E\Prob(E)\E[M\ind_R\mid E]\\
 &=\frac1P\sum_E\Prob(E)\E[M\mid E]\Prob(R\mid E)\\
 &=\sum_Ew_Ep_E\\
 &\ge\frac14\sum_E\frac{w_E}{d_E}.
\end{align*}
By Cauchy--Schwarz
\begin{align*}
 q^2
 &=\left(\sum_E\sqrt{\frac{w_E}{d_E}}\sqrt{w_Ed_E}\right)^2\\
 &\le\left(\sum_E\frac{w_E}{d_E}\right)\left(\sum_Ew_Ed_E\right)\\
 &=G\sum_E\frac{w_E}{d_E},
\end{align*}
and thus
\[
 \frac{\E[M\ind_R]}P\ge\frac{q^2}{4G}.
\]

It remains to estimate the gain of $\BS$ on good samples. To do this we will fix a sample $\calX=\{X_1,\dots,X_n\}$, assume it is good, and randomize only its arrival times. Since the sample is good, the highest occupied bucket contains a unique observation, which has value $M$. Assign independent uniform arrival times in $[0,1]$ to the observations in $\calX$, and recall that $\BS$ rejects before time $a$ and then accepts the first observation whose bucket is strictly above every previously seen bucket.

Condition on $M$ arriving at time $u\ge a$, and on the set $A$ of observations that arrive before time $u$. First consider the case where $A\neq\emptyset$, and let $r$ be the number of observations in the highest bucket occupied by an observation in $A$. The arrival times of these $r$ observations are independent and uniform on $[0,u]$, and $\BS$ selects $M$ if and only if at least one of them arrives before time $a$. Thus
\begin{align*}
 \Prob(\text{$\BS$ selects $M$}\mid u, \text{$\calX$ good}, A) 
 =1-\left(1-\frac au\right)^r
 \ge\frac au.
\end{align*}
If $A=\emptyset$, the left-hand side is equal to~$1$ and the inequality also holds.

Suppose now that $\calX$ belongs to $R$, and let $x$ and $y$ be two arbitrary observations in the second-highest occupied bucket. Conditional on $M$ arriving at time $u$,
\[
 \Prob(\text{$x$ and $y$ arrive before time $u$}\mid u, \text{$\calX$ in $R$})=u^2.
\]
Whenever both $x$ and $y$ arrive before time $u$, the highest bucket occupied by an element of $A$ is the second-highest occupied bucket for the entire sample, so $r\ge2$ and thus
\begin{align*}
\Prob(\text{$\BS$ selects $M$}\mid u, \text{$\calX$ in $R$}, A, \text{$x$ and $y$ arrive before time $u$}) 
&= 1-\left(1-\frac au\right)^r \\
&\ge 1-\left(1-\frac au\right)^2
=\frac au+\frac au\left(1-\frac au\right).
\end{align*}
By averaging over the possible choices of $A$,
\begin{align*}
 \Prob(\text{$\BS$ selects $M$}\mid u,\text{$\calX$ good})
 \ge\frac au+\ind_R\,u^2\frac au\left(1-\frac au\right).
\end{align*}
By integrating over the arrival time $u$ of $M$,
\begin{align*}
 \Prob(\text{$\BS$ selects $M$}\mid\text{$\calX$ good})
 &\ge \int_a^1\frac au\,du
       +\ind_R\int_a^1a(u-a)\,du\\
 &= -a\log a+\ind_R\frac{a(1-a)^2}{2}\\
 &= a+\kappa\ind_R.
\end{align*}
Finally, by multiplying by $M$ when $\BS$ selects $M$, and averaging over good samples,
\begin{align*}
 \frac{\E[\BS]}{P}
 &\ge a\frac{\E[M\ind_{\rm good}]}P
       +\kappa\frac{\E[M\ind_R]}P\\
 &\ge aq+\frac{\kappa q^2}{4G}.
\end{align*}
Subtracting $a$ from both sides yields the claim.
\end{proof}

\subsection{Large Gap Helps Threshold Secretary}
\label{sec:threshold}

Draw $T$ uniformly at random from the $N$ cuts, and let 
\[
 z=G/N,\qquad b=ae^{-z}.
\]
Since $G\le Nq\le N$, it follows that $z\in[0,1]$ and $b\in[a^2,a]$. The threshold secretary rule $\TS$ rejects before time $b$, accepts the first record above $T$ in $[b,a)$, and accepts the first record from time $a$ onward. If $G=0$, then $b=a$ and $\TS$ behaves like the secretary rule. When the gap $G$ is large, a random cut will often isolate the maximum value from all other values, allowing $\TS$ to select it before time~$a$. 
We obtain the following improvement over the secretary rule.
\begin{lemma}
\label{lem:threshold}
For every instance,
\begin{equation*}
 \frac{\E[\TS]}{P}-a\ge a(e^{-z}+z-1)+a^2e^{-z}\gamma.
\end{equation*}
\end{lemma}
\begin{proof}
Fix a sample and a cut~$T$. 
If $T$ isolates the maximum from all other values, $\TS$ selects $M$ whenever $\SEC_a$ does, but in addition selects $M$ with certainty if it arrives in $[b,a)$; it thus selects $M$ with probability at least $2a-b$.
If $T$ does not isolate the maximum, $\TS$ still selects $M$ with probability at least $-b\log b$: indeed, it selects $M$ whenever $\SEC_b$ does if $T\leq M$, and behaves like $\SEC_a$ if $T>M$.
Exactly $D$ of the $N$ cuts isolate the maximum on a good sample, and none do on a bad sample, so 
\[
 \frac{\E\!\left[
 M\ind_{\{\text{cut isolates maximum}\}}
 \right]}{P}
 =
 \frac{\E\!\left[M\ind_{\rm good}\,D/N\right]}{P}
 =
 \frac{G}{N}
 =
 z.
\]
Moreover, for every sample and cut, every value selected by $\SEC_b$ after time~$a$ is also selected by $\TS$, so $\TS$ selects the second maximum with probability at least
\[
 b\int_a^1\frac{1-u}{u}\,du=ab.
\]
Thus
\begin{align*}
 \frac{\E[\TS]}{P}-a
 &\geq z(2a-b) + (1-z)(-b\log b) + ab\gamma - a \\
 &= z(2a-b+b\log b) - b\log b + ab\gamma - a \\
 &\geq z(a-b) - b\log b + ab\gamma - a \\
 & = za(1-e^{-z}) + ae^{-z}(1+z) + a^2e^{-z}\gamma - a \\
  &= a(e^{-z}+z-1)+a^2e^{-z}\gamma,
\end{align*}
where the first inequality holds by taking expectation over possible samples and cuts and because selection of the maximum and second maximum are disjoint events, the second inequality because $-b\log b\leq a$ and thus $b\log b\geq -a$, the second equality because $b=ae^{-z}$ and $\log a=-1$, and the other two equalities by rearranging. This shows the claim.
\end{proof}

\subsection{Bucket Secretary or Threshold Secretary Gains}
\label{sec:gain}

To achieve the lower bound in \Cref{thm:main}, the gambler determines~$G$ and uses either $\BS$ or $\TS$ depending on which is better. By \Cref{lem:bucket}, $\BS$ improves over the secretary rule's guarantee of $1/e$ by order $G^{-1}$. By \Cref{lem:threshold}, $\TS$ improves over it by order $(G/N)^2$. The larger of the two improvements is minimized when they are equal, which happens when $G\approx N^{2/3}$ for a guarantee of order $N^{-2/3}$. The additional logarithm enters through the number of cuts per component distribution, which is logarithmic.
\begin{lemma}[Lower bound]
\label{lem:lower}
There exists an absolute constant $c>0$ such that, for all $k\geq 2$,
\[
 c_k \ge \frac1e+\frac{c}{(k\log k)^{2/3}}.
\]
\end{lemma}
\begin{proof}
For $n\le2$, choosing a uniformly random observation achieves a competitive ratio of at least $1/2$ and the claim holds. Now assume $n\ge3$, let $L=\lceil\log_2 k\rceil$
and
\[
 r=\frac{\max\{\E[\BS],\E[\TS]\}}{P}-a.
\]
Then
\[
 \delta\le\frac1k,
 \qquad
 N\le k(\lceil\log_2 k\rceil+1)=O(k\log k).
\]
As before let $z=G/N$, so that $0\le z\le1$ and thus
\begin{align}
\label{eq:az2}
e^{-z}+z-1
 &=\int_0^z(1-e^{-t})\,dt \notag\\
 &=\int_0^z\int_0^t e^{-s}\,ds\,dt \notag\\
 &=\int_0^z\int_s^z e^{-s}\,dt\,ds \notag\\
 &=\int_0^z(z-s)e^{-s}\,ds \notag\\
 &=\int_0^1 z(1-t)e^{-zt}\,z\,dt \notag\\
 &=z^2\int_0^1(1-t)e^{-zt}\,dt \geq az^2, 
\end{align}
where for the third equality we have reversed the order of integration over
$0\le s\le t\le z$, the fifth equality holds by substituting $s=zt$ if $z>0$ and because both sides are zero if $z=0$, and the inequality holds because the integral decreases in $z$ and equals $a$ at $z=1$.
It follows that
\begin{equation}
\label{eq:threshold}
\begin{aligned}
 r &\geq a(e^{-z}+z-1)+a^2e^{-z}\gamma \\
 &\geq a^2z^2+a^3\gamma,
\end{aligned}
\end{equation}
where the first inequality holds by \Cref{lem:threshold}, and the second inequality by~\Cref{eq:az2} and because $e^{-z}\ge a$.

If $\gamma\ge1/4$ then
\[
 r\ge a^3/4,
\]
and the claim holds.

Now consider the case where $\gamma\leq 1/4$. Then, for $k\ge4$,
\[
 q\ge 1-\gamma-\delta\ge 1-\frac{1}{4} - \frac{1}{k} \geq \frac{1}{2}
\]
where the first inequality holds by \Cref{lem:loss}. In particular $G\ge q>0$, and thus
\[
 r \geq \frac{\kappa q^2}{4G}-a(1-q)
 \geq \frac{\kappa}{16G}-a(\gamma+\delta),
\]
where the first inequality holds by \Cref{lem:bucket}.
Multiplying this inequality by~$a^2$ and adding~\Cref{eq:threshold} yields
\begin{align}
 \label{eq:bucket}
 (a^2+1)r
 &\ge
 a^2\left(\frac{\kappa}{16G}-a(\gamma+\delta)\right) + a^2z^2+a^3\gamma \notag\\
 &\ge
 a^2\left(\frac{\kappa}{16G}+z^2\right)
 -a^3\delta \notag\\
 &= a^2\left(\frac{\kappa}{16G}+\frac{G^2}{N^2}\right)
 -a^3\delta.
\end{align}
For fixed $N>0$ and any $G>0$ let 
\[
  f(G) = \frac{\kappa}{16G}+\frac{G^2}{N^2}.
\]
Then
\[
 f'(G)=-\frac{\kappa}{16G^2}+\frac{2G}{N^2},
 \qquad
 f''(G)=\frac{\kappa}{8G^3}+\frac{2}{N^2}>0,
\]
so $f$ has a unique global minimum where
\[
 \frac{\kappa}{16G}=\frac{2G^2}{N^2},
\]
and thus 
\[
 G = \left(\frac{\kappa N^2}{32}\right)^{1/3}
\]
This value minimizes the right-hand side of~\Cref{eq:bucket}, so
\begin{align*}
 (a^2+1)r 
 &\geq a^2\left( \frac{\kappa}{16(\kappa N^2/32)^{1/3}} +\frac{(\kappa N^2/32)^{2/3}}{N^2} \right)-a^3\delta \\
 &= a^2\kappa^{2/3}N^{-2/3} \left(\frac{32^{1/3}}{16}+\frac1{32^{2/3}}\right)-a^3\delta \\
 &= a^2\kappa^{2/3}N^{-2/3} \left(\frac2{32^{2/3}}+\frac1{32^{2/3}}\right)-a^3\delta \\
 &= 3a^2\left(\frac{\kappa}{32}\right)^{2/3}N^{-2/3}-a^3\delta.
\end{align*}
Thus, for absolute constants $c_1,c_2>0$,
\[
 r\geq c_1(k\log k)^{-2/3}-\frac{c_2}{k}.
\]
Since $1/k=o\bigl((k\log k)^{-2/3}\bigr)$, setting $c=c_1/2$ proves the claim when $k\geq k_0$ for some absolute constant $k_0\geq 4$. To see that it also holds when $k<k_0$, we note that a mixture of $k$ components can be turned into a mixture of $k_0$ components by adding copies of one of the components, so that $c_k\geq c_{k_0}>1/e$. Since there are only finitely many values $k<k_0$, the claim thus holds for all $k\geq 2$ for a smaller value of~$c$.
\end{proof}

\section{Upper Bound}
\label{sec:upper}

We prove the upper bound using a mixture in which the latent component determines only the scale of the rewards. To make the largest reward dominate all others, we place the possible rewards at powers of a base $b>1$, and call $m$ the \emph{level} of the reward $b^m$; we will eventually let $b$ tend to infinity. For a parameter $\lambda>0$ chosen later, let $\rho=e^{-\lambda/n}$,
and let $G_1,\dots,G_n$ be independent of one another and of the latent
component~$J$, with
\[
  \mathbb{P}(G_i=g)=(1-\rho)\rho^g,\qquad g=0,1,\dots.
\]
Now let component distribution $F_j$ be the distribution of $b^{j-G_i}$, so that the observed rewards are $X_i=b^{J-G_i}$. Call the largest possible reward $b^j$ in component~$j$ the \emph{endpoint} of that component. The true endpoint is then a random variable $Z=b^J$, and $G_i$ counts the number of levels by which draw~$i$ falls below $Z$. Each draw is the endpoint with probability $1-\rho$, so the endpoint appears at least once in $n$ draws with probability $1-\rho^n=1-e^{-\lambda}$, independently of~$n$. For any integer $m\le j$,
\begin{equation}\label{eq:levels}
  \mathbb{P}(X_i=b^m\mid J=j)=(1-\rho)\rho^{j-m},\qquad
  \mathbb{P}(X_i\le b^m\mid J=j)=\rho^{j-m},
\end{equation}
where the second equality holds because $\mathbb{P}(G_i\ge j-m)=\rho^{j-m}$. Moving from $F_j$ to $F_{j+1}$ thus multiplies the probability of every old level by the same factor~$\rho$. This limits what the gambler can learn about the component from low observations. Each observation $X_i=Z\,b^{-G_i}$ is the product of the endpoint and a factor $b^{-G_i}\le1$, where endpoint and factors are independent. Our mixture is thus an instance of the common-scale model of \citet{CorreaEtAl2026ThresholdDynamics}, with the additional property that the factors are identically distributed and the scale takes exactly~$k$ values.

When observing a record, a value that exceeds all previous values, the gambler must decide between selecting that record and waiting for a larger value. For this decision to be difficult, it should matter across many components, which can be achieved by giving components with a larger endpoint a smaller prior probability. Consider weights $w(1),\dots,w(k)>0$ such that $\sum_{j=1}^kw(j)=1$, to be chosen later. We set the prior probabilities so that the conditional probability of selecting the endpoint of component \(j\) receives relative weight \(w(j)\) in the expected payoff. If component \(j\) has prior probability \(\pi_j\), then selecting its endpoint with conditional probability \(q_j\) contributes \(\pi_j b^j q_j\) to the expected payoff. We thus define 
\[
  \pi_j=\frac{w(j)b^{-j}}{\sum_{\ell=1}^k w(\ell)b^{-\ell}}.
\]

In \Cref{sec:endpoint} we reduce optimal stopping on this mixture to optimal stopping in an auxiliary problem where the gambler observes levels and seeks to maximize the probability of selecting the level of the endpoint. Rather than solving the auxiliary problem, in \Cref{sec:surplus} we bound its value by a quantity we call the stopping surplus, averaged over time and summed over levels. The stopping surplus measures how much better it is to stop than to continue when observing level $m$ as a record at time $s\in[0,1]$, the fraction of draws observed so far. The reduction and the stopping surplus bound apply for any choice of the weights $w(1),\dots,w(k)$. In \Cref{sec:weights} we identify weights for which the summed surplus is at most $1/e$ plus two error terms, of order $\lambda$
and $1/(k\lambda)^2$. In \Cref{sec:balancing} we finally choose $\lambda$ to balance those terms and obtain the desired upper bound.

\subsection{Reduction to Endpoint Selection}
\label{sec:endpoint}

The weights $w(1),\dots,w(k)$ have a natural interpretation in the context of an auxiliary stopping problem where the gambler observes levels and seeks to maximize the probability of selecting the endpoint level. 
Let $J$ be a random variable such that $\Prob(J=\ell)=w(\ell)$, and let $Y_1,\dots,Y_n$ be drawn independently according to the distribution over levels in component $F_J$. The gambler observes $Y_1,\dots,Y_n$ sequentially and selects one of them at the time it is revealed, with the goal of maximizing the probability of selecting value $J$, \ie the level of the endpoint of $F_J$. The best such probability is
\[
 V_n=\sup_\tau\Prob(\tau\le n,\ Y_\tau=J),
\]
where the supremum is taken over random variables~$\tau$ such that the event $\{\tau=i\}$ for $i=1,\dots,n$ is determined by the observations $Y_1,\dots,Y_i$ and any private randomness of the gambler. The following lemma establishes that the prophet problem for the mixture is as difficult as the auxiliary endpoint selection problem, scaled by the probability $1-e^{-\lambda}$ with which the endpoint appears.
\begin{lemma}\label{lem:upper-endpoint}
For every $n$, every choice of weights $w(1),\dots,w(k)>0$ with $\sum_{j=1}^kw(j)=1$, and $\lambda>0$,
\begin{equation*}%
  c_k\le \frac{V_n}{1-e^{-\lambda}}.
\end{equation*}
\end{lemma}
\begin{proof}
Fix $b>1$, let $J$ be drawn from prior $\pi$ and $Z=b^J$. Since $b>1$, observing $X_i=b^{J-G_i}$ is equivalent to observing $J-G_i$, so the mixture and the auxiliary problem admit the same stopping rules and conditional on $J=\ell$ the distribution over levels is the same in both. Moreover $\E[Z]=\sum_{j=1}^k\pi_jb^j=1/\sum_{j=1}^kw(j)b^{-j}$, and thus $\pi_j b^j=w(j)\,\E[Z]$ for $1\le j\le k$.
For any stopping rule $\tau$,
\begin{align*}
  \E[X_\tau]
  &= \sum_{j=1}^k\pi_j\,\E[X_\tau\mid J=j] \\
  & \le\sum_{j=1}^k\pi_j b^j\Bigl(\Prob(\tau\le n,\,G_\tau=0\mid J=j)+\frac1b\Bigr)\\
  &= \E[Z]\Bigl(\sum_{j=1}^kw(j)\,\Prob(\tau\le n,\,G_\tau=0\mid J=j)+\frac1b\Bigr)\\
  &\le \E[Z]\Bigl(V_n+\frac1b\Bigr),
\end{align*}
where the first inequality holds because, for $J=j$, $X_\tau=b^j$ if $G_\tau=0$ and $X_\tau\leq b^{j-1}$ if $G_\tau\ge 1$, the second equality because $\pi_j b^j=w(j)\,\E[Z]$ and $\sum_{j=1}^kw(j)=1$, and the second inequality because the event $G_\tau=0$ corresponds to selection of the endpoint.
Also
\[
  \E[M]=\sum_{j=1}^k\pi_j\,\E[M\mid J=j]
  \ge\sum_{j=1}^k\pi_j b^j\bigl(1-\rho^n\bigr)
  =\bigl(1-e^{-\lambda}\bigr)\,\E[Z],
\]
where the inequality holds because, for $J=j$, $M=b^j$ if $G_i=0$ for
some $i\le n$, which happens with probability $1-\rho^n$, and $M\ge 0$
otherwise, and the second equality because $\pi_j b^j=w(j)\,\E[Z]$, $\sum_{j=1}^kw(j)=1$, and
$\rho^n=e^{-\lambda}$.
Since $M\le Z\le b^k$, it follows that $0<\E[M]<\infty$, so the mixture is one of those over which the infimum in the definition of $c_k$ is taken. Thus
\[
  c_k\le\frac{\sup_\tau\E[X_\tau]}{\E[M]}
  \le\frac{V_n+1/b}{1-e^{-\lambda}},
\]
The auxiliary problem, and thus $V_n$, does not depend on~$b$, so letting $b$ go to infinity shows the claim.
\end{proof}

\subsection{The Value of an Opportunity To Stop}
\label{sec:surplus}

We now analyze the auxiliary problem with endpoint probabilities $w(1),\dots,w(k)$. We measure time by the fraction of draws observed so far, so that draw~$i$ occurs at time $s=i/n\in[0,1]$. Rather than solving this problem optimally, we bound the value of being allowed to stop with level $m$ at time $s$ from above by comparing it to the value of waiting for the next level that is strictly higher than $m$.

The comparison matters only when $m$ may be the endpoint, \ie if it is at least~$1$ and at least as high as every previously observed level. In that case, assume that the true endpoint is $m+\ell$ for $\ell\ge 0$. The probability that none of the $\ell$ levels above~$m$ has appeared by time~$s$ is then $e^{-\lambda s\ell}$. Since the levels below and including $m$ appear with the same relative probabilities for every value of $\ell$, the posterior probabilities of the endpoints $m,\dots,m+\ell$ are proportional to $w(m),\dots,w(m+\ell)e^{-\ell\lambda s}$ and stopping with level $m$ at time $s$ succeeds with probability proportional to $w(m)$. Waiting to instead select the first level above~$m$ succeeds if such a level appears after time~$s$, which happens with probability $1-e^{-\lambda(1-s)\ell}$, and if it is the endpoint, which happens with probability at least $1/\ell$ because the endpoint is the most likely of the $\ell$ levels above~$m$. Multiplying these probabilities and summing over~$\ell$ shows that waiting succeeds with at least $L_m(s)/w(m)$ times the probability of stopping, where
\begin{equation}\label{eq:Lm}
 L_m(s)=\sum_{\ell\geq 1}w(m+\ell)\frac{e^{-\lambda s\ell}-e^{-\lambda \ell}}{\ell}
\end{equation}
with $w(j)=0$ for $j>k$.

Having the option to stop helps the gambler only if stopping beats waiting, and then by at most the success probability of stopping times
$1-L_m(s)/w(m)$. Call
\[
  \bigl(w(m)-L_m(s)\bigr)_+=w(m)\Bigl(1-\frac{L_m(s)}{w(m)}\Bigr)_+
\]
the \emph{stopping surplus} at level~$m$ and time~$s$, where $x_+=\max\{x,0\}$. Since $L_m(s)$ bounds the value of waiting from below, the stopping surplus bounds how much the option to stop can add to the gambler's probability of success. The following lemma bounds $V_n/\lambda$, the success probability per expected appearance of the endpoint in the auxiliary problem, by the stopping surplus summed over levels and averaged over time.
\begin{lemma}\label{lem:upper-policy}
For every $k$, every choice of weights $w(1),\dots,w(k)>0$ with $\sum_{j=1}^kw(j)=1$, and $\lambda>0$,
\[
 \limsup_{n\to\infty}\frac{V_n}{\lambda}
 \leq \sum_{m=1}^k\int_0^1\bigl(w(m)-L_m(s)\bigr)_+\dd s.
\]
\end{lemma}
\begin{proof}
Write $Y_i=J-G_i$ for the level of draw~$i$. By~\Cref{eq:levels}, $\Prob(Y_i=y\mid J=j)=(1-\rho)\rho^{j-y}$ and $\Prob(Y_i\le y\mid J=j)=\rho^{j-y}$ for every integer $y\le j$. 
Let $p_i=\Prob(J=Y_i\mid Y_1,\dots,Y_i)$ be the probability that draw~$i$ is the endpoint, $r_i$ the optimal probability, conditional on $Y_1,\dots,Y_i$, of selecting the endpoint among draws $i+1,\dots,n$, and $v_i$ the optimal probability of selecting the endpoint when all draws are observed but only draws $i,\dots,n$ may be selected. Then $v_1=V_n$ and $v_{n+1}=0$. Moreover $v_i=\E[\max\{p_i,r_i\}]$ and $v_{i+1}=\E[r_i]$, since a gambler who may select draw~$i$ can stop or continue optimally, and a gambler who may not select it can only continue. Thus
\begin{equation}\label{eq:Vn}
  V_n=\sum_{i=1}^n(v_i-v_{i+1})=\sum_{i=1}^n\E\bigl[(p_i-r_i)_+\bigr].
\end{equation}

Stopping at draw~$i$ can succeed only if $Y_i=\max\{1,Y_1,\dots,Y_i\}$. Thus $p_i=0$ outside the disjoint events
\[
  E_{i,m}=\{Y_i=m,\ Y_1,\dots,Y_{i-1}\le m\}, \qquad 1\leq m\leq k.
\]
Fix $i$, $m$, and levels $y_1,\dots,y_i$ with $y_i=m$ and
$y_1,\dots,y_{i-1}\le m$. Then $\Prob(Y_1=y_1,\dots,Y_i=y_i\mid J=j)=0$ for $j<m$, and for
$0\le\ell\le k-m$,
\[
  \Prob(Y_1=y_1,\dots,Y_i=y_i\mid J=m+\ell)=\prod_{u=1}^i(1-\rho)\rho^{m+\ell-y_u}
  =\rho^{i\ell}\prod_{u=1}^i(1-\rho)\rho^{m-y_u}.
\]
The last product does not depend on~$\ell$, so for all such levels $y_1,\dots,y_i$, the posterior probability of $J=m+\ell$ given $Y_1=y_1,\dots,Y_i=y_i$ is proportional to $w(m+\ell)\rho^{i\ell}$. Conditional on $J$, the remaining draws are independent of $Y_1,\dots,Y_i$, so $p_i$ and $r_i$ are constant on $E_{i,m}$. Moreover, by~\Cref{eq:levels}, $\Prob(J=m+\ell,\,E_{i,m})=w(m+\ell)(1-\rho)\rho^{i\ell}$.

To bound $\E[r_i\mathbf{1}_{E_{i,m}}]$ from below, consider the stopping rule that selects the first
draw after~$i$ whose level is strictly greater than~$m$. Conditional on $J=m+\ell$ for $\ell\ge1$, it follows from~\Cref{eq:levels} and
independence that such a draw exists with probability $1-\rho^{\ell(n-i)}$ and is the endpoint with probability $(1-\rho)/(1-\rho^\ell)$. We claim that
\begin{align}
  \label{eq:piri}
  \E\bigl[(p_i-r_i)_+\mathbf{1}_{E_{i,m}}\bigr]
  &=\bigl(\E[p_i\mathbf{1}_{E_{i,m}}]-\E[r_i\mathbf{1}_{E_{i,m}}]\bigr)_+ \notag\\
  &=\bigl(\Prob(J=m,\,E_{i,m})-\E[r_i\mathbf{1}_{E_{i,m}}]\bigr)_+ \notag\\
  &\le\Bigl(\Prob(J=m,\,E_{i,m})-\sum_{\ell=1}^{k-m}\Prob(J=m+\ell,\,E_{i,m})
     \bigl(1-\rho^{\ell(n-i)}\bigr)\frac{1-\rho}{1-\rho^\ell}\Bigr)_+ \notag\\
  &\le(1-\rho)\Bigl(w(m)-\sum_{\ell=1}^{k-m}w(m+\ell)\frac{\rho^{i\ell}-\rho^{n\ell}}{\ell}\Bigr)_+ \notag\\
  &=(1-\rho)\bigl(w(m)-L_m(i/n)\bigr)_+,
\end{align}
Indeed, the first equality holds because $p_i-r_i$ is constant on $E_{i,m}$. The second equality holds because $Y_i=m$ whenever $E_{i,m}$ occurs, and thus $\E[p_i\mathbf{1}_{E_{i,m}}]=\Prob(J=Y_i,\,E_{i,m})=\Prob(J=m,\,E_{i,m})$. 
The first inequality holds because conditional on $Y_1,\dots,Y_i$, the rule succeeds with probability at most $r_i$, so that $\E[r_i\mathbf{1}_{E_{i,m}}]$ is at least the probability that $E_{i,m}$ occurs and the rule succeeds.
The second inequality holds because $\Prob(J=m+\ell,\,E_{i,m})=w(m+\ell)(1-\rho)\rho^{i\ell}$ and
$(1-\rho)/(1-\rho^\ell)=1/(1+\rho+\dots+\rho^{\ell-1})\ge1/\ell$. The third equality finally holds by~\Cref{eq:Lm}, and because $\rho^{i\ell}=e^{-\lambda(i/n)\ell}$, $\rho^{n\ell}=e^{-\lambda \ell}$, and $w(j)=0$ for $j>k$. Now
\begin{align*}
  V_n &= \sum_{i=1}^n\sum_{m=1}^k\E\bigl[(p_i-r_i)_+\mathbf{1}_{E_{i,m}}\bigr] \\
  &\le n(1-\rho)\sum_{m=1}^k\frac1n\sum_{i=1}^n\bigl(w(m)-L_m(i/n)\bigr)_+ \\
  &\le \lambda\sum_{m=1}^k\frac1n\sum_{i=1}^n\bigl(w(m)-L_m(i/n)\bigr)_+,
\end{align*}
where the equality holds by~\Cref{eq:Vn}, because $p_i=0$ outside
the events $E_{i,m}$, and $r_i\ge0$, the first inequality by~\Cref{eq:piri}, and the second inequality because
$n(1-\rho)=n(1-e^{-\lambda/n})\le\lambda$. $L_m$ is a finite sum of continuous functions, so for each~$m$ the
inner average is a Riemann sum of a continuous function on $[0,1]$ and
converges to its integral. This proves the claim.
\end{proof}

\subsection{Choosing Weights}
\label{sec:weights}

\Cref{lem:upper-policy} bounds the gambler's success probability for any choice of the reward weights. We now choose the weights so that the integral in this bound is close to $1/e$. Recall that $L_m(s)$ is a weighted sum over levels above $m$. To estimate it, we will approximate $w(m+\ell)$ by $w(m)+\ell\,w'(m)$. A symmetric choice of weights will let us pair levels of equal weight and opposite slope. We also want this approximation to remain accurate near the last component, where the weights at higher levels are zero.

Let 
\begin{equation*}
w(x)= \begin{cases}
\displaystyle\frac{x^2(k+1-x)^2}{Z_k}, & 0\le x\le k+1,\\[4pt]
0, & \text{otherwise},
\end{cases}
\end{equation*}
where $Z_k=\sum_{j=1}^k j^2(k+1-j)^2$. We will refer to these weights, which have degree four in~$x$ as the \emph{quartic weights}. The polynomial $x(k+1-x)$ is symmetric around $(k+1)/2$ and vanishes at $0$ and $k+1$. Squaring it makes the slope vanish at both ends as well, so it can be extended by zero without a jump in slope. Division by $Z_k$ implies that $\sum_{j=1}^k w(j)=1$.

The following lemma lets us bound the error when approximating $w(m+\ell)$ by $w(m)+\ell\,w'(m)$.
\begin{lemma}\label{lem:upper-weights}
For the quartic weights $w$, $1\le m\le k$, and all integers $\ell\ge1$,
\[
 \bigl|w(m+\ell)-w(m)-\ell w'(m)\bigr|=O(\ell^2/k^3).
\]
Moreover \(w(k+1-m)=w(m)\), \(w'(k+1-m)=-w'(m)\), and
\begin{equation*}
 \sum_{m=1}^k\frac{w'(m)^2}{w(m)}
   =\frac{4}{Z_k}\sum_{m=1}^k(k+1-2m)^2=O(k^{-2}).
\end{equation*}
\end{lemma}
\begin{proof}
First observe that
\begin{equation}
\label{eq:zbound}
Z_k=\sum_{j=1}^k j^2(k+1-j)^2
\ge \frac{k}{2}\left(\frac{k+1}{4}\right)^4
\ge \frac{k^5}{512},
\end{equation}
where the first inequality holds because there are at least \(k/2\) integers \(j\in[(k+1)/4,\,3(k+1)/4]\), and for each of them both \(j\) and \(k+1-j\) are at least \((k+1)/4\).
By differentiation, for \(0<x<k+1\),
\[
 w'(x)=\frac{2x(k+1-x)(k+1-2x)}{Z_k},\qquad
 w''(x)=\frac{12(x-(k+1)/2)^2-(k+1)^2}{Z_k}.
\]
The numerator of \(w''(x)\) has absolute value at most \(2(k+1)^2\), so~\Cref{eq:zbound} yields
\begin{equation}
  \label{eq:curvature}
  \lvert w''(x)\rvert\le \frac{C}{k^3}    
\end{equation}
for an absolute constant \(C>0\). For \(x<0\) or \(x>k+1\), \(w(x)=0\) and thus \(w'(x)=0\).

Fix \(m\in\{1,\ldots,k\}\) and an integer \(\ell\ge 1\). For \(0\le t\le \ell\), let \(u_t=\min\{m+t,k+1\}\). Then
\begin{align*}
\bigl|w(m+\ell)-w(m)-\ell w'(m)\bigr|
&=\left|\int_0^\ell \bigl(w'(m+t)-w'(m)\bigr)\,dt\right|\\
&\le \int_0^\ell \bigl|w'(u_t)-w'(m)\bigr|\,dt\\
&\le \int_0^\ell\int_m^{u_t}|w''(x)|\,dx\,dt\\
&\le \frac{C}{k^3}\int_0^\ell (u_t-m)\,dt\\
&\le \frac{C}{k^3}\int_0^\ell t\,dt
= \frac{C\ell^2}{2k^3},
\end{align*}
where the first equality holds because $w$ is continuously differentiable and thus \(w(m+\ell)-w(m)=\int_0^\ell w'(m+t)\,dt\), the first inequality by the triangle inequality and because either \(u_t=m+t\) or \(u_t=k+1<m+t\) and thus $w'(u_t)=w'(m+t)=0$, the second inequality by the triangle inequality and because $w'$ agrees with the derivative of the polynomial on $[m,u_t]$ and thus \( w'(u_t)-w'(m)=\int_m^{u_t} w''(x)\,dx\), the third inequality by~\Cref{eq:curvature}, and the fourth inequality because \(u_t\le m+t\) and thus \(u_t-m\le t\). This proves the first claim for all \(\ell\ge 1\), including those for which \(m+\ell>k+1\).

That \(w(k+1-x)=w(x)\) follows directly from the definition of $w$, and differentiating yields \(-w'(k+1-x)=w'(x)\).

Finally
\begin{align*}
\sum_{m=1}^k\frac{w'(m)^2}{w(m)}
&= \sum_{m=1}^k \frac{4m^2(k+1-m)^2(k+1-2m)^2}{Z_k^2}
\cdot\frac{Z_k}{m^2(k+1-m)^2} \\
&=\frac{4}{Z_k}\sum_{m=1}^k(k+1-2m)^2 \\
&\le \frac{4k(k-1)^2}{Z_k} \\
&\le \frac{2048}{k^2},
\end{align*}
where the first inequality holds because \(\lvert k+1-2m\rvert\le k-1\) and the second inequality by~\Cref{eq:zbound}.
\end{proof}

We are now ready to bound the integrated stopping surplus. For an indication of the leading term, suppose that in~\Cref{eq:Lm} $w(m+\ell)$ was replaced by $w(m)$. At time $s>0$, the sum of the coefficients would then equal \(\log((1-e^{-\lambda})/(1-e^{-\lambda s}))\), which tends to \(-\log s\) as \(\lambda\) tends to zero. The corresponding surplus would approach \(w(m)(1+\log s)_+\). This limiting surplus is positive precisely after the secretary cutoff \(1/e\), and its integral is \(w(m)/e\). To bound how far the actual integral can exceed this value, we must account for positive \(\lambda\) and for the variation of the weights. Pairing equal weights on opposite sides of the midpoint will cancel their slopes except near the cutoff. We obtain the following lemma.
\begin{lemma}\label{lem:upper-integral}
For the quartic weights $w$ and \(0<\lambda\le1\),
\[
 \sum_{m=1}^k\int_0^1(w(m)-L_m(s))_+\,ds
 \le \frac1e+O\!\left(\lambda+\frac{1}{k^2\lambda^2}\right).
\]
\end{lemma}
\begin{proof}
We first bound $\bigl(w(m)-L_m(s)\bigr)_+$ for each level $m$ and time~$s$, then sum the bounds while pairing levels $m$ and $k+1-m$, and finally integrate over~$s$.
  
Fix $s\in[0,1]$. By \Cref{lem:upper-weights}, there exists an absolute constant $C>0$ such that, for $1\le m\le k$ and integers $\ell\ge 1$,
\begin{equation}\label{eq:lem7}
  \bigl|w(m+\ell)-w(m)-\ell w'(m)\bigr|\le\frac{C\ell^2}{k^3}
  \qquad\text{and}\qquad
  \sum_{m=1}^k\frac{w'(m)^2}{w(m)}\le\frac{C}{k^2}.
\end{equation}
Let $\sbar=\max\{s,e^{-2}\}$. For $\ell\ge 1$, let
\[
  a_\ell=\frac{e^{-\lambda \sbar \ell}-e^{-\lambda \ell}}{\ell}.
\]
Let
\[
  A_0=\sum_{\ell\ge 1}a_\ell
  \qquad\text{and}\qquad
  A_1=\sum_{\ell\ge1}\ell\,a_\ell.
\]
Then
\begin{align}
  \label{eq:A0bound}
  A_0 &= \sum_{\ell\ge1}\frac{e^{-\lambda \sbar \ell}}{\ell}-\sum_{\ell\ge1}\frac{e^{-\lambda \ell}}{\ell} \notag\\
  &= \log\bigl(1-e^{-\lambda}\bigr)-\log\bigl(1-e^{-\lambda\sbar}\bigr) \notag\\
  &\ge \log\bigl(\lambda e^{-\lambda}\bigr)-\log(\lambda\sbar) \notag\\
  &= -\log \sbar-\lambda,
\end{align}
where the second equality holds because $\sum_{\ell\ge1}z^\ell/\ell=-\log(1-z)$ when $0\le z<1$, and the inequality because $xe^{-x}\le1-e^{-x}\le x$ when $x\ge 0$. Also
\[
  A_1=\sum_{\ell\ge 1}\bigl(e^{-\lambda \sbar \ell}-e^{-\lambda \ell}\bigr),
\]
and thus
\begin{equation}
  \label{eq:A1bounds}
  0\le A_1\le\sum_{\ell\ge1}e^{-\lambda \sbar \ell}=\frac{1}{e^{\lambda \sbar}-1}
  \le\frac{1}{\lambda \sbar}\le\frac{e^2}{\lambda},
\end{equation}
where the first inequality holds because $a_\ell\ge 0$, the third inequality because $e^x-1\ge x$ when $x\ge0$, and
the fourth inequality because $\sbar\ge e^{-2}$. Finally,
\begin{align}
  \label{eq:sumi1bound}
  \sum_{\ell\ge1}\ell^2a_\ell
  \le\sum_{\ell\ge1}\ell e^{-\lambda \sbar \ell}
  =\frac{e^{-\lambda \sbar}}{(1-e^{-\lambda \sbar})^2}
  =\frac{1}{(e^{\lambda \sbar/2}-e^{-\lambda \sbar/2})^2}
  \le\frac{1}{(\lambda \sbar)^2}
  \le\frac{e^4}{\lambda^2},
\end{align}
where the first inequality holds because
$\ell^2a_\ell=\ell(e^{-\lambda \sbar \ell}-e^{-\lambda \ell})\le \ell e^{-\lambda \sbar \ell}$, the first
equality because $\sum_{\ell\ge1}\ell z^\ell=z/(1-z)^2$ when $0\le z<1$, the second
inequality because $e^x-e^{-x}\ge 2x$ when $x\ge0$, and the third inequality
because $\sbar\ge e^{-2}$.

Now fix $m$ with $1\le m\le k$. For $\ell\ge1$, let $R_\ell=w(m+\ell)-w(m)-\ell w'(m)$. Then
\begin{align}
  \label{eq:pointwise}
  \bigl(w(m)-L_m(s)\bigr)_+
  &\le\bigl(w(m)-L_m(\sbar)\bigr)_+\notag\\
  &=\Bigl(w(m)\bigl(1-A_0\bigr)-w'(m)A_1-\sum_{\ell\ge1}R_\ell\,a_\ell\Bigr)_+\notag\\
  &\le\Bigl(w(m)(1+\log \sbar+\lambda)-w'(m)A_1+\frac{C}{k^3}\sum_{\ell\ge1}\ell^2a_\ell\Bigr)_+\notag\\
  &\le\Bigl(w(m)(1+\log \sbar)-w'(m)A_1+w(m)\lambda+\frac{Ce^4}{k^3\lambda^2}\Bigr)_+\notag\\
  &\le\bigl(w(m)(1+\log \sbar)-w'(m)A_1\bigr)_++w(m)\lambda+\frac{Ce^4}{k^3\lambda^2},
\end{align}
where the first inequality holds because every summand in~\Cref{eq:Lm} is nonincreasing in~$s$ and $\sbar\ge s$, so that $L_m(\sbar)\le L_m(s)$, the second inequality because $1-A_0\le 1+\log \sbar+\lambda$ by~\Cref{eq:A0bound} and because $-R_\ell\le |R_\ell|\le C\ell^2/k^3$ by~\Cref{eq:lem7}, the third inequality by~\Cref{eq:sumi1bound}, and the fourth inequality because $(x+y)_+\le x_++y$ when $y\ge 0$.

In summing over $m$ we will use that, for all $x$ and $y$,
\begin{equation}\label{eq:plusplus}
  (x-y)_++(x+y)_+=2x_++\bigl(|y|-|x|\bigr)_+.
\end{equation}
Indeed, both sides are unchanged when $y$ is replaced by $-y$, so we may
assume that $y\ge0$. If $y\le|x|$, then $x-y$ and $x+y$ are both nonnegative
when $x\ge0$ and both nonpositive when $x<0$, so both sides equal $2x_+$. If
$y>|x|$, then $x-y<0<x+y$, so both sides equal $x+y$.

For $1\le m\le k$, let $x_m=w(m)(1+\log \sbar)$ and $y_m=w'(m)A_1$. Then
\begin{align}
  \label{eq:sumbound}
  \sum_{m=1}^k(x_m-y_m)_+
  &=\frac12\sum_{m=1}^k\Bigl((x_m-y_m)_++(x_{k+1-m}-y_{k+1-m})_+\Bigr)\notag\\
  &=\frac12\sum_{m=1}^k\Bigl((x_m-y_m)_++(x_m+y_m)_+\Bigr)\notag\\
  &=\sum_{m=1}^k(x_m)_++\frac12\sum_{m=1}^k\bigl(|y_m|-|x_m|\bigr)_+\notag\\
  &=(1+\log \sbar)_++\frac12\sum_{m=1}^k\bigl(|w'(m)|A_1-w(m)|1+\log \sbar|\bigr)_+,
\end{align}
where the first equality holds because $m\mapsto k+1-m$ is a permutation of
$\{1,\dots,k\}$, the second equality because, by \Cref{lem:upper-weights}, $w(k+1-m)=w(m)$ and
$w'(k+1-m)=-w'(m)$, the third equality by~\Cref{eq:plusplus}, and the fourth equality because $w(m)>0$ when
$1\le m\le k$, $\sum_{m=1}^kw(m)=1$, and $A_1\ge 0$. 
The $m$th summand in the second sum is positive precisely when $|w'(m)|A_1>w(m)|1+\log \sbar|$.

We now integrate~\Cref{eq:sumbound} over $s\in[0,1]$. For the first term on the right-hand side,
\begin{equation}
  \label{eq:term1}
  \int_0^1(1+\log \sbar)_+\,ds=\int_{1/e}^1(1+\log s)\,ds
  =\bigl[s\log s\bigr]_{1/e}^1=\frac1e,
\end{equation}
where the first equality holds because $\sbar\le 1/e$, and thus $(1+\log \sbar)_+=0$ when $s\le 1/e$, and because $\sbar=s$ for $s\ge 1/e$.
To bound the second term, we use that $|1+\log \sbar|\ge|s-1/e|$ for all $s\in[0,1]$. For $s\le e^{-2}$ this holds because $|1+\log \sbar|=1>|s-1/e|$. For $s\ge e^{-2}$ it holds because $\sbar=s$, because $1+\log(1/e)=0$, and because
the derivative of $1+\log s$ is $1/s\ge1$. Thus, for $1\le m\le k$,
\begin{align}
  \label{eq:term2}
  \int_0^1\bigl(|w'(m)|A_1-w(m)|1+\log \sbar|\bigr)_+\,ds
  &\le \int_0^1\Bigl(\frac{e^2|w'(m)|}{\lambda}-w(m)\Bigl|s-\frac1e\Bigr|\Bigr)_+\,ds \notag\\
  &\le \int_{-\infty}^{\infty}\Bigl(\frac{e^2|w'(m)|}{\lambda}-w(m)\Bigl|s-\frac1e\Bigr|\Bigr)_+\,ds \notag\\
  &= \frac{e^4w'(m)^2}{\lambda^2w(m)},
\end{align}
where the first inequality holds because $A_1\le e^2/\lambda$ and $|1+\log \sbar|\ge|s-1/e|$, the second inequality because the integrand is nonnegative, and the equality because the integrand is a triangle centered at $1/e$ with height $h=e^2|w'(m)|/\lambda$ and base $2h/w(m)$.

We finally obtain
\begin{align*}
  \sum_{m=1}^k\int_0^1\bigl(w(m)-L_m(s)\bigr)_+\,ds
  &\le\int_0^1\sum_{m=1}^k\bigl(w(m)(1+\log \sbar)-w'(m)A_1\bigr)_+\,ds
     +\lambda+\frac{Ce^4}{k^2\lambda^2}\\
  &=\int_0^1(1+\log \sbar)_+\,ds
     +\frac12\sum_{m=1}^k\int_0^1\bigl(|w'(m)|A_1-w(m)|1+\log \sbar|\bigr)_+\,ds\\
  &\qquad+\lambda+\frac{Ce^4}{k^2\lambda^2}\\
  &\le\frac1e+\frac{e^4}{2\lambda^2}\sum_{m=1}^k\frac{w'(m)^2}{w(m)}
     +\lambda+\frac{Ce^4}{k^2\lambda^2}\\
  &\le\frac1e+\lambda+\frac{3Ce^4}{2k^2\lambda^2},
\end{align*}
where the first inequality holds by~\Cref{eq:pointwise}, because $\sum_{m=1}^kw(m)=1$, and because the last term in~\Cref{eq:pointwise} does not depend on $m$, the equality by~\Cref{eq:sumbound}, the second inequality by~\Cref{eq:term1} and~\Cref{eq:term2}, and the third inequality by~\Cref{eq:lem7}. This proves the claim.
\end{proof}

\subsection{Balancing Excess}
\label{sec:balancing}

We will now prove the upper bound in \Cref{thm:main} by balancing two effects that are both controlled by $\lambda$ and contribute to the excess $c_k-1/e$.
On the one hand, a smaller value of $\lambda$ makes the endpoint less likely to appear. This hurts the gambler because, conditional on the endpoint appearing at all, the gambler is less likely to get a second chance to select it. The excess associated with positive endpoint frequency is of order $\lambda$. On the other hand, a smaller $\lambda$ also means that more distant levels contribute to the continuation value. Near the secretary cutoff $s=1/e$, the probability \(e^{-\lambda s \ell}\) that all draws so far have been at or below $m$ when the latent endpoint is $m+\ell$ decays over distances $\ell$ of order $1/\lambda$. The weights vary on the scale of $k$ levels, so the relevant distance is a fraction $1/(k\lambda)$ of that scale. Pairing levels with equal weights and opposite slopes cancels the first-order variation up to a residual whose integral is quadratic, which means the excess is of order $1/(k\lambda)^2$. 
Balancing the two error bounds amounts to choosing $\lambda$ such that $\lambda=(k\lambda)^{-2}$, which gives $\lambda=k^{-2/3}$.
\begin{lemma}[Upper bound]
\label{lem:upper}
There exists an absolute constant $C>0$ such that, for all $k\geq 1$,
\[
  c_k \leq \frac{1}{e}+\frac{C}{k^{2/3}}.
\]
\end{lemma}
\begin{proof}
Use the quartic weights and $\lambda=k^{-2/3}$. As before, let $V_n$ denote the maximum probability with which the gambler can select the endpoint when there are $n$ draws. Then
\begin{align*}
\frac{1-e^{-\lambda}}{\lambda}\,c_k
&\le \limsup_{n\to\infty}\frac{V_n}{\lambda} \\
&\le \sum_{m=1}^k\int_0^1\bigl(w(m)-L_m(s)\bigr)_+\,ds \\
&\le \frac1e+C_0\left(\lambda+\frac{1}{k^2\lambda^2}\right) \\
&= \frac1e+2C_0\lambda
\end{align*}
for an absolute constant \(C_0>0\). Indeed, the first inequality holds because $c_k(1-e^{-\lambda})\le V_n$ for every $n$ by \Cref{lem:upper-endpoint}, 
the second inequality by \Cref{lem:upper-policy}, the third inequality by \Cref{lem:upper-integral}, and the equality because $\lambda=k^{-2/3}$.
Thus
\begin{align*}
  c_k &\le \frac{\lambda}{1-e^{-\lambda}} \left(\frac1e+2C_0\lambda\right) \\
  &\le (1+\lambda)\left(\frac1e+2C_0\lambda\right) \\
  &\le \frac1e+\left(\frac1e+4C_0\right)\lambda \\
  &= \frac1e+\left(\frac1e+4C_0\right) k^{-2/3},
\end{align*}
where the second inequality holds because $e^\lambda\ge 1+\lambda$ and thus $1-e^{-\lambda}\ge \lambda/(1+\lambda)$, the third inequality because $\lambda=k^{-2/3}\le 1$ for all $k\geq 1$, and the equality because $\lambda=k^{-2/3}$. The claim follows for $C=1/e+4C_0$.
\end{proof}

\paragraph{Acknowledgments}
The first author was supported by the Engineering and Physical Sciences Research Council under grant UKRI2406 ``Prophet Inequalities for Mixtures.'' Large language models Claude and ChatGPT were used as tools in the discovery of the results and the production of the writeup. Full responsibility for the work lies with the authors.

\bibliographystyle{abbrvnat}

\end{document}